\documentclass[11pt]{article}
\usepackage{graphicx}
\usepackage{amsfonts}
\usepackage{theorem}
\newtheorem{Lem}{Lemma}[section]

\newtheorem{The}[Lem]{Theorem}
\newtheorem{Prop}[Lem]{Proposition}
\newtheorem{Cor}[Lem]{Corollary}

\newtheorem{Rem}[Lem]{Remark}

\newcommand{\qed}{\hbox{\rule{6pt}{6pt}}}

\begin{document}
\title{Refined Young inequalities with Specht's ratio}
\author{Shigeru Furuichi$^1$\footnote{E-mail:furuichi@chs.nihon-u.ac.jp}  \\
$^1${\small Department of Computer Science and System Analysis,}\\
{\small College of Humanities and Sciences, Nihon University,}\\
{\small 3-25-40, Sakurajyousui, Setagaya-ku, Tokyo, 156-8550, Japan}}
\date{}
\maketitle

{\bf Abstract.} In this paper, we show that the $\nu$-weighted arithmetic mean is greater than the product of the $\nu$-weighted geometric mean and Specht's ratio.
As a corollary, we also show that the $\nu$-weighted geometric mean is greater than the product of the $\nu$-weighted harmonic mean and Specht's ratio.
These results give the improvements for the classical Young inequalities, since Specht's ratio is generally greater than $1$.
In addition, we give an operator inequality for positive operators, applying our refined Young inequality.
\vspace{3mm}


\vspace{3mm}
\vspace{3mm}
\maketitle

\section{Introduction}
We start from the famous Young inequality:
\begin{equation}  \label{young_ineq}
(1-\nu) a + \nu b \geq a^{1-\nu} b^{\nu}
\end{equation}
for positive numbers $a$, $b$ and $\nu \in [0,1]$. The inequality (\ref{young_ineq})  is also called $\nu$-weighted  arithmetic-geometric mean inequality
and its reverse inequality was given in \cite{Tom1} with Specht's ratio as follows:
\begin{equation}  \label{rev_young_ineq}
S\left(\frac{a}{b}\right) a^{1-\nu} b^{\nu} \geq (1-\nu) a + \nu b 
\end{equation}
for positive numbers $a$, $b$ and $\nu \in [0,1]$. 
Where the Specht's ratio \cite{Specht,JIFUJII} was defined by
$$
S(h) \equiv \frac{h^{\frac{1}{h-1}}}{e\log{h^{\frac{1}{h-1}}}},\quad (h \neq 1)
$$
for a positive real number $h$. 
%Note that $S(h) > 1$ for $h>0, h\neq 1$ and $S(1)=1$.

Recently, based on the refined Young inequality \cite{BK,KM}:
\begin{equation} \label{improve_ineq}
(1-\nu) a+ \nu b \geq a^{1-\nu} b^{\nu}+ r (\sqrt{a}-\sqrt{b})^2,
\end{equation}
for positive numbers $a$, $b$ and $\nu \in [0,1]$, where $r\equiv \min\{\nu,1-\nu\}$,
we proved the following operator inequalities:
\begin{Prop}  \label{Prev_Furuichi} {\bf (\cite{Furuichi})}
For $\nu \in [0,1]$ and positive operators $A$ and $B$, we have
\begin{eqnarray*}  
 (1-\nu )A+\nu B  &\geq& A\sharp_{\nu} B  +  2r \left(\frac{A+B}{2} -A\sharp_{1/2}B \right) \label{the01-ineq01} \\
&\geq & A\sharp_{\nu} B  \label{the01-ineq02} \\
&\geq& \left\{ A^{-1}\sharp_{\nu} B^{-1}+2r \left(\frac{A^{-1}+B^{-1}}{2} -A^{-1}\sharp_{1/2}B^{-1} \right)\right\}^{-1} \label{the01-ineq03} \\
&\geq& \left\{  (1-\nu )A^{-1}+\nu B^{-1}   \right\}^{-1}   \label{the01-ineq04} 
\end{eqnarray*}
where $r\equiv \min\left\{\nu,1-\nu\right\}$ and  $A\sharp_{\nu}B \equiv A^{1/2}(A^{-1/2}BA^{-1/2})^{\nu}A^{1/2}$  defined for $\nu\in[0,1]$.
\end{Prop}
The above inequalities can be regarded as an additive-type refinement for the Young inequalities \cite{FY,Furuta}:
\begin{equation}  \label{young_ineq_op}
(1-\nu )A+\nu B  \geq A\sharp_{\nu} B  \geq \left\{  (1-\nu )A^{-1}+\nu B^{-1}   \right\}^{-1}.   
\end{equation}

In this short paper, we give a multiplicative-type refinement for the Young inequalities (\ref{young_ineq_op}) with the Specht's ratio.

%%%%%%%%%%%%%%%%%%%%%%%%%%%%%%%%%%%%%%%%%%%%%%%%%%%%%%%%%%%%%%%%%%%%%%%%%%%%%%%%%%%%%%%%%%%%%%%%%%%%%%%%%%%%%%%%%%%%%%%%%%%%%%%%%%%%%%%%%%%%%%%%%%%%%%%%%%%%%%%
%%%%%%%%%%%%%%%%%%%%%%%%%%%%%%%%%%%%%%%%%%%%%%%%%%%%%%%%%%%%%%%%%%%%%%%%%%%%%%%%%%%%%%%%%%%%%%%%%%%%%%%%%%%%%%%%%%%%%%%%%%%%%%%%%%%%%%%%%%%%%%%%%%%%%%%%%%%%%%%
%%%%%%%%%%%%%%%%%%%%%%%%%%%%%%%%%%%%%%%%%%%%%%%%%%%%%%%%%%%%%%%%%%%%%%%%%%%%%%%%%%%%%%%%%%%%%%%%%%%%%%%%%%%%%%%%%%%%%%%%%%%%%%%%%%%%%%%%%%%%%%%%%%%%%%%%%%%%%%%
%%%%%%%%%%%%%%%%%%%%%%%%%%%%%%%%%%%%%%%%%%%  Section2  %%%%%%%%%%%%%%%%%%%%%%%%%%%%%%%%%%%%%%%%%%%%%%%%%%%%%%%%%%%%%%%%%%%%%%%%%%%%%%%%%%%%%%%%%%%%%%%%%%%%%%%%
%%%%%%%%%%%%%%%%%%%%%%%%%%%%%%%%%%%%%%%%%%%%%%%%%%%%%%%%%%%%%%%%%%%%%%%%%%%%%%%%%%%%%%%%%%%%%%%%%%%%%%%%%%%%%%%%%%%%%%%%%%%%%%%%%%%%%%%%%%%%%%%%%%%%%%%%%%%%%%%
%%%%%%%%%%%%%%%%%%%%%%%%%%%%%%%%%%%%%%%%%%%%%%%%%%%%%%%%%%%%%%%%%%%%%%%%%%%%%%%%%%%%%%%%%%%%%%%%%%%%%%%%%%%%%%%%%%%%%%%%%%%%%%%%%%%%%%%%%%%%%%%%%%%%%%%%%%%%%%%
%%%%%%%%%%%%%%%%%%%%%%%%%%%%%%%%%%%%%%%%%%%%%%%%%%%%%%%%%%%%%%%%%%%%%%%%%%%%%%%%%%%%%%%%%%%%%%%%%%%%%%%%%%%%%%%%%%%%%%%%%%%%%%%%%%%%%%%%%%%%%%%%%%%%%%%%%%%%%%%
\section{Main results}
We here review the properties of the Specht's ratio. See \cite{Tom1,Specht,JIFUJII} for example, as for the proof and the details.
\begin{Lem}  \label{s-p}
The Specht's ratio 
$$S(h) \equiv \frac{h^{\frac{1}{h-1}}}{e\log{h^{\frac{1}{h-1}}}},\quad (h \neq 1,\,\, h>0)$$
has the following properties. 
\begin{itemize}
\item[(i)] $S(1) = 1$ and $S(h) = S(1/h) >1 $ for $h>0$.
\item[(ii)] $S(h)$ is a monotone increasing function on $(1,\infty)$.
\item[(iii)] $S(h)$ is a monotone decreasing function on $(0,1)$.
\end{itemize}
\end{Lem}

We use the following lemmas to show our theorem.
\begin{Lem} \label{lem01}
For $x\geq 1$, we have
\begin{equation}  \label{ineq02}
\frac{2(x-1)}{x+1} \leq \log x \leq \frac{x-1}{\sqrt{x}}.
\end{equation}
%For $0< x \leq 1$, we have
%\begin{equation}  \label{ineq02a}
% \frac{x-1}{\sqrt{x}} \leq \log x \leq \frac{2(x-1)}{x+1}.
%\end{equation}
\end{Lem}

{\it Proof}:
We firstly prove the second inequality of (\ref{ineq02}).
We put $\sqrt{x}=t$ and 
$$
f(t) \equiv \frac{t^2-1}{t}-2\log t, \quad (t \geq 1).
$$
Then we have $f'(t)=\left(\frac{t-1}{t}\right)^2\geq 0$ and $f(1)=0$.
Thus we have $f(t) \geq f(1) =0$ and then we have $\log t^2 \leq  \frac{t^2-1}{t}$,
which implies the second inequality in (\ref{ineq02}).

We also put
$$
g(x) \equiv (x+1) \log x -2(x-1),\quad (x\geq 1).
$$
Then we have $f(1)=0$, $f'(x)=\log x+\frac{x+1}{x}-2$, $f'(1)=0$ and $f''(x)=\frac{t-1}{t^2}\geq 0$.
Therefore we have $f(x) \geq f(1) =0$, which implies the first inequality in (\ref{ineq02}).
%The inequalities (\ref{ineq02a}) are proven by putting $x=\frac{1}{y}$ in the inequalities  (\ref{ineq02}).

\hfill \qed

Note that Lemma \ref{lem01} can be also proven by the following relation for three means:
$$
\sqrt{xy} < \frac{x-y}{\log x-\log y} < \frac{x+y}{2}
$$
for positive real numbers $x$ and $y$, where $x \neq y$.

\begin{Lem} \label{lem-e-t}
For $t>0$, we have
\begin{equation} \label{lemma-e-t}
e (t^2+1) \geq (t+1) t^{\frac{t}{t-1}}.
\end{equation}
\end{Lem}
{\it Proof}:
We firstly prove the inequality (\ref{lemma-e-t}) for $t \geq 1$.
We put
$$
f(t) = e (t^2+1) - (t+1) t^{\frac{t}{t-1}}.
$$
By using the first inequality of (\ref{ineq02}), we have 
\begin{eqnarray*}
f'(t)&=& \frac{2t(t-1)^2e+2t(1-t)t^{\frac{t}{t-1}} +t^{\frac{t}{t-1}}(t+1)\log t }{(t-1)^2}\\
&\geq & \frac{2t(t-1)^2e+2t(1-t)t^{\frac{t}{t-1}} +2(t-1)t^{\frac{t}{t-1}} }{(t-1)^2}\\
&=& \frac{2t(t-1)^2e- 2t(t-1)^2t^{\frac{1}{t-1}}  }{(t-1)^2}\\
&\geq & \frac{2t(t-1)^2t^{\frac{1}{t-1}}   - 2t(t-1)^2  t^{\frac{1}{t-1}}  }{(t-1)^2}\\
&=&0.
\end{eqnarray*}
In the last inequality, we have used the fact that $\lim_{t\to 1}t^{\frac{1}{t-1}}=e$ and the function $t^{\frac{1}{t-1}}$ is monotone decreasing
 on $t \in [1,\infty)$.
We also have $f(1)=0$ so that we have $f(t) \geq 0$ which proves the following inequality:
$$ e (t^2+1) \geq (t+1) t^{\frac{t}{t-1}},\quad  t \geq 1.$$
Putting $t=\frac{1}{s}$ in the above inequality with simple calculations, we have
$$ e (s^2+1) \geq (s+1) s^{\frac{s}{s-1}},\quad  0 < s \leq 1.$$

\hfill \qed

Then we have the following inequality which improves the classical Young inequality 
between $\nu$-weighted geometric mean and $\nu$-weighted arithmetic mean.

\begin{The} \label{the-r-y}
For $a,b>0$ and $\nu\in[0,1]$,
\begin{equation} \label{improve_ineq_mult_scalar}
(1-\nu) a + \nu b \geq S\left(\left(\frac{b}{a}\right)^r\right) a^{1-\nu} b^{\nu},
\end{equation}
where $r \equiv \min \left\{\nu,1-\nu\right\}$ ans $S(\cdot)$ is the Specht's ratio.
\end{The}

{\it Proof}:
We prove the following inequality
\begin{equation}
\frac{(b-1)\nu +1}{b^{\nu}S(b^{\nu})} = \frac{e\left\{(b-1)\nu +1\right\}\log b^{\nu}}{\left(b^{\nu}\right)^{\frac{b^{\nu}}{b^{\nu}-1}}(b^{\nu}-1)} \geq 1
\end{equation}
 in the case of $0 \leq \nu \leq \frac{1}{2}$.
From Lemma \ref{lem01}, we have
$$
\frac{\log b^{\nu}}{b^{\nu}-1} \geq \frac{2}{b^{\nu} +1}, \quad b>0 .
$$
Therefore we have the following first inequality:
\begin{equation}
 \frac{e\left\{(b-1)\nu +1\right\}\log b^{\nu}}{\left(b^{\nu}\right)^{\frac{b^{\nu}}{b^{\nu}-1}}(b^{\nu}-1)} 
\geq \frac{2e\left\{(b-1)\nu +1\right\}}{\left(b^{\nu}\right)^{\frac{b^{\nu}}{b^{\nu}-1}}(b^{\nu}+1)} \geq 1,
\end{equation}
thus we have only to prove the above second inequality.
For this purpose, we put the following function $f_b$ on  $\nu\in[0,\frac{1}{2}]$ for  $b>0$:
$$
f_b(\nu) \equiv 2e \left\{ (b-1)\nu +1\right\} -\left(b^{\nu}\right)^{\frac{b^{\nu}}{b^{\nu}-1}}(b^{\nu}+1).
$$
Then we have
\begin{eqnarray*}
&& f''_b(\nu)=-\frac{\left(\log b\right)^2}{(b^{\nu}-1)^4}\left(b^{\nu}\right)^{\frac{2b^{\nu}-1}{b^{\nu}-1}} \left\{(b^{\nu}-1)^2(4b^{2\nu}-5b^{\nu}-1) \right.\\
&& \left.\hspace*{25mm} -(b^{\nu}-1)^2(3b^{\nu}+1)\log b^{\nu} + b^{\nu} (b^{\nu}+1)(\log b^{\nu})^2\right\}.
\end{eqnarray*}
For the case of $b \geq 1$, using the inequalities (\ref{ineq02}), we have
\begin{eqnarray*}
&& (b^{\nu}-1)^2(4b^{2\nu}-5b^{\nu}-1)  -(b^{\nu}-1)^2(3b^{\nu}+1)\log b^{\nu} + b^{\nu} (b^{\nu}+1)(\log b^{\nu})^2 \\
&& \geq (b^{\nu}-1)^2(4b^{2\nu}-5b^{\nu}-1) -(b^{\nu}-1)^2(3b^{\nu}+1) \frac{b^{\nu}-1}{b^{\nu/2}} + b^{\nu} (b^{\nu}+1) 
\left( \frac{2(b^{\nu}-1)}{b^{\nu}+1}\right)^2\\
&&=\frac{   (b^{\nu/2}-1)^4 (b^{\nu/2}+1)^3 (4b^{2\nu}+b^{3\nu/2}+4 b^{\nu} +1 )  }{b^{\nu/2}    (b^{\nu}+1)} \geq 0
\end{eqnarray*}
For the case of $0< b \leq 1$, using the inequalities  (\ref{ineq02}), we also have
\begin{eqnarray*}
&& (b^{\nu}-1)^2(4b^{2\nu}-5b^{\nu}-1)  -(b^{\nu}-1)^2(3b^{\nu}+1)\log b^{\nu} + b^{\nu} (b^{\nu}+1)(\log b^{\nu})^2 \\
&& = (b^{\nu}-1)^2(4b^{2\nu}-5b^{\nu}-1)  +(b^{\nu}-1)^2(3b^{\nu}+1)\log \frac{1}{b^{\nu}} + b^{\nu} (b^{\nu}+1)\left(\log \frac{1}{b^{\nu}}\right)^2 \\
&& \geq (b^{\nu}-1)^2(4b^{2\nu}-5b^{\nu}-1) +(b^{\nu}-1)^2(3b^{\nu}+1)\left( \frac{2(\frac{1}{b^{\nu}}-1)}{\frac{1}{b^{\nu}}+1}\right) + b^{\nu} (b^{\nu}+1) 
\left( \frac{2  \left(  \frac{1}{b^{\nu}}-1 \right) }{ \frac{1}{b^{\nu}} +1 }  \right) ^2\\
&&=\frac{(b^{\nu}-1)^4(4b^{\nu}+1) }{b^{\nu}+1}  \geq 0
\end{eqnarray*}
Thus we have $f''_b(\nu) \leq 0$ for $b>0$.
In addition, we have $f_b(0)=0$ and $f_b(\frac{1}{2})=e(b+1)-(\sqrt{b}+1)(\sqrt{b})^{\frac{\sqrt{b}}{\sqrt{b}-1}} \geq 0$,
applying Lemma \ref{lem-e-t} with $t=\sqrt{b} >0$. Therefore we have $f_b(\nu) \geq 0$ for $\nu \in [0,\frac{1}{2}]$.
Thus we have the following inequality
\begin{equation} \label{the-final01}
\frac{(b-1)\nu +1}{b^{\nu}S(b^{\nu})} \geq 1,\quad 0\leq \nu \leq \frac{1}{2},\,\, b>0
\end{equation}
which implies
$$
 \nu b + (1-\nu)  \geq S(b^{\nu}) b^{\nu}.
$$
Replacing $b$ by $\frac{b}{a}$ in the above inequality and then multiplying $a$ to the both sides, we have
$$
 (1-\nu)a + \nu b \geq S\left(\left(\frac{b}{a}\right)^{\nu}\right)a^{1-\nu}b^{\nu},\quad 0 \leq \nu \leq \frac{1}{2}, \,\,a,b>0
$$

Finally, from the inequality (\ref{the-final01}), we have
$$
\frac{(a-1)\mu +1}{a^{\mu}S(a^{\mu})} \geq 1,\quad 0\leq \mu \leq \frac{1}{2}, \,\, a>0.
$$
Putting $\nu = 1- \mu$ in the above inequality we have
$$
\nu +(1-\nu)a \geq a^{1-\nu} S(a^{1-\nu}), \quad \frac{1}{2} \leq \nu \leq 1, \,\,b>0.
$$
Replacing $a$ by $\frac{a}{b}$ in the above inequality and then multiplying $b$ to the both sides, we have
$$
 (1-\nu) a+ \nu b \geq S\left(\left(\frac{a}{b}\right)^{1-\nu}\right)a^{1-\nu}b^{\nu},\quad \frac{1}{2}  \leq \nu \leq 1,\,\, a,b>0,
$$
since $S(1/h) = S(h)$ for $h>0$, ((i) of Lemma \ref{s-p}).
Thus the proof of the present theorem was completed.

\hfill \qed

\begin{Rem}
Theorem \ref{the-r-y} gives a tighter lower bound of the $\nu$-wighted arithmetic mean of two variables, 
since the Specht's ratio is greater than $1$, ((i) of Lemma \ref{s-p}).
\end{Rem}

The following inequality also improves the relation between $\nu$-weighted geometric mean and $\nu$-weighted harmonic mean.

\begin{Cor}  \label{Cor}
For positive numbers $a$, $b$ and $\nu \in[0,1]$, we have 
\begin{equation}\label{ineq04}
S\left(\left(\frac{a}{b}\right)^r\right) \left( (1- \nu) \frac{1}{a} +  \nu \frac{1}{b} \right)^{-1}\leq a^{1-\nu}b^{\nu},
\end{equation} 
where $r\equiv \min\left\{\nu,1-\nu \right\}$ and  $S(\cdot)$ is the Specht's ratio.
\end{Cor}
{\it Proof}:
Replace $a$ and $b$ in Theorem \ref{the-r-y} by $\frac{1}{a}$ and $\frac{1}{b}$, respectively.

\hfill \qed

Applying Theorem \ref{the-r-y}, we have the following operator inequality for positive operators.
\begin{The} \label{the-r-y-op}
For two positive operators $A$, $B$ and positive real numbrs $m,m',M,M'$ satisfying the following conditions (i) or (ii):
\begin{itemize}
\item[(i)]$0<m'I \leq A \leq mI<MI\leq B \leq M'I$ 
\item[(ii)]$0<m'I \leq B \leq mI<MI\leq A \leq M'I$ 
\end{itemize}
with $h\equiv \frac{M}{m}$ and $h'\equiv\frac{M'}{m'}$, we have
\begin{eqnarray}
(1-\nu) A + \nu B &\geq& S\left(h^r\right) A\sharp_{\nu}B \label{the-01}  \\
&\geq& A\sharp_{\nu}B \label{the-02}  \\
&\geq& S\left(h^r\right) \left\{(1-\nu)A^{-1} +\nu B^{-1} \right\}^{-1}  \label{the-03}  \\
 &\geq& \left\{(1-\nu)A^{-1} +\nu B^{-1} \right\}^{-1},\label{the-04}  
\end{eqnarray}
where  $\nu\in [0,1]$,  $r\equiv \min \left\{\nu,1-\nu\right\}$, $S(\cdot)$ is the Specht's ratio  and 
$A\sharp_{\nu}B \equiv A^{1/2}\left(A^{-1/2}BA^{-1/2}\right)^{\nu}A^{1/2}$is the $\nu$-power mean for positive operators $A$ and $B$ \cite{KA}.
\end{The}

{\it Proof}:
From Theorem \ref{the-r-y}, we have
$$
\nu x+(1-\nu) \geq S(x^r) x^{\nu}
$$
for any $x>0$.
Therefore we have
$$
\nu X+(1-\nu)I \geq \min_{m'\leq x \leq M'} S(x^r) X^{\nu}
$$
for the positive operator $X$ such that $0<m'I \leq X \leq M'I$.
We here put $X=A^{-1/2}BA^{-1/2}$.

In the case of (i),  we have $h=\frac{M}{m}\leq A^{-1/2}BA^{-1/2} \leq \frac{M'}{m'} =h'$.
Then we have
$$
\nu A^{-1/2}BA^{-1/2}+ (1-\nu)I \geq \min_{h\leq x \leq h'} S(x^r) \left(A^{-1/2}BA^{-1/2}   \right)^{\nu}.
$$
Since $S(x)$ is an increasing function for $x >1$, ((ii) of Lemma \ref{s-p})
we have
\begin{equation} \label{same_ineq}
\nu A^{-1/2}BA^{-1/2}+ (1-\nu)I  \geq S(h^r) \left(A^{-1/2}BA^{-1/2}   \right)^{\nu}.
\end{equation}

In the case of (ii),
we also have $\frac{1}{h'}=\frac{m'}{M'} \leq A^{-1/2}BA^{-1/2} \leq \frac{m}{M}=\frac{1}{h}$.
Then we also have
$$
\nu A^{-1/2}BA^{-1/2}+ (1-\nu)I \geq \min_{\frac{1}{h'} \leq x \leq \frac{1}{h}} S(x^r) \left(A^{-1/2}BA^{-1/2}   \right)^{\nu}.
$$
Since $S(x)$ is a decreasing function for $0<x<1$ ((iii) of Lemma \ref{s-p}),
we have 
$$
\nu A^{-1/2}BA^{-1/2}+ (1-\nu)I  \geq S\left(\frac{1}{h^r}\right) \left(A^{-1/2}BA^{-1/2}   \right)^{\nu}.
$$
By the property $S(x)=S(1/x)$ for $x>0$ ((i) of Lemma \ref{s-p}),  the above inequality is the same to (\ref{same_ineq}).
Multiplying $A^{1/2}$ from the both sides to the inequality (\ref{same_ineq}), we have the inequality (\ref{the-01}).

The inequality (\ref{the-03}) can be proven by replacing $A$ and $B$ by $A^{-1}$ and $B^{-1}$, respectively in the
first inequality and taking its inverse.

The inequality (\ref{the-02}) and the inequality (\ref{the-04}) are trivial, due to the property of the Specht's ratio
$S(x) \geq 1$ for $x >0$.

\hfill \qed

\section{Conclusion}
We have shown the refined Young inequalities for a real number with Specht ratio.
Applying these inequalities we have obtained their operator version inequalities which refine 
the classical Young operator inequalities as our previous results have done in Proposition \ref{Prev_Furuichi} (See \cite{Furuichi}).
Therefore we have two different refinements for  the classical Young inequalities (\ref{young_ineq_op}). 
Two kinds of the operator inequalities are based on the scalar inequalities (\ref{improve_ineq}) and (\ref{improve_ineq_mult_scalar}).

In our previous paper \cite{Furuichi}, we have proved the additive-type refined Young inequality for $n$ real numbers.
\begin{Prop}  \label{the_gen01}  {\bf (\cite{Furuichi})}
Let  $a_1,\cdots,a_n\geq 0$ and $p_1,\cdots,p_n > 0$ with $\sum_{j=1}^n p_j=1$ and  $\lambda  \equiv \min \left\{ p_1,\cdots ,p_n \right\}$.
If we assume that the multiplicity attaining $\lambda$ is $1$,
then we have
\begin{equation} \label{gen_ineq01}
\sum_{i=1}^n p_i a_i  -   \prod_{i=1}^n a_i^{p_i}  \geq n \lambda  \left(\frac{1}{n}\sum_{i=1}^n a_i-  \prod_{i=1}^n a_i^{1/n}     \right),
\end{equation}
with equality if and only if $a_1=\cdots =a_n$.
\end{Prop}
See \cite{Mit,Ald} for recent developments based on the above inequality (or Jensen-type inequality \cite{Dra}). 
It is also notable that we do not need the assumption that the multiplicity attaining $\lambda$ is $1$, to prove only inequality (\ref{gen_ineq01}).
This assumption connects with the equality condition.

Closing this section, we give comments on the multiplicative-type refined Young inequality for $n$ real numbers.
We have not yet found its proof. We also have not found any counter-examples for the following  3-variables case:
\begin{equation} \label{multiplicative_type_ineq02}
w_1 a_1 +w_2 a_2 +w_3 a_3 \geq S(h^r) a_1^{w_1} a_2^{w_2} a_3^{w_3},
\end{equation}
for $a_i \in [m,M]$ where $0<m<M$ with $h\equiv \frac{\max\left\{a_1,a_2,a_3\right\}}{\min\left\{a_1,a_2,a_3\right\}}$ 
and $r \equiv \min\left\{ w_1,w_2,w_3 \right\}$, where $w_i > 0$ and $w_1+w_2+w_3=1$.

The problem on the multiplicative-type refined Young inequality for $n$ real numbers will be our future work.


\begin{thebibliography}{99}
\bibitem{Tom1} M. Tominaga, Specht's ratio in the Young inequality, Sci.Math.Japon.,Vol.55(2002),pp.583-588.
\bibitem{Specht} W.Specht, Zer Theorie der elementaren Mittel, Math.Z.,Vol.74(1960),pp.91-98.
\bibitem{JIFUJII} J.I.Fujii,S.Izumino and Y.Seo, Determinant for positive operators and Specht's theorem, Sci.Math.Japon.,Vol.1(1998),pp.307-310.
\bibitem{BK} N.A.Bobylev and M. A. Krasnoselsky, Extremum Analysis (degenerate cases), Moscow, preprint, 1981, 52 pages, (in Russian).
\bibitem{KM} F.Kittaneh and Y.Manasrah, Improved Young and Heinz inequalities for matrices, J.Math.Anal.Appl.,Vol.36(2010), pp.262-269.
\bibitem{Furuichi} S.Furuichi, On refined Young inequalities and reverse inequalities, J.Math.Ineq..,Vol.5(2011), pp.21-31.
\bibitem{FY} T.Furuta and M.Yanagida, Generalized means and convexity of inversion for positive operators, Amer.Math.Monthly, Vol.105 (1998),pp.258-259.
\bibitem{Furuta} T.Furuta, Invitation to linear operators: From matrix to bounded linear operators on a Hilbert space, Taylor and  Francis, 2002.
\bibitem{KA} F.Kubo and T.Ando, Means of positive operators, Math. Ann.,Vol.264(1980),pp.205-224.
\bibitem{Mit}F. C. Mitori, About the precision in Jensen-Steffensen inequality,  Annals of the University of Craiova, Mathematics and Computer Science Series, Vol.37 (2010), pp.73-84.
\bibitem{Ald}J. M. Aldaz, Comparison of differences between arithmetic and geometric means, arXiv:1001.5055v2.
\bibitem{Dra} S.S.Dragomir, Bounds for the normalised Jensen functional, Bull. Australian Math. Soc., Vol.74(2006), pp.471-478.    
\end{thebibliography}
\end{document}